\documentclass{amsart}
\usepackage{setspace}
\usepackage{a4}
\usepackage{amssymb,amsmath,amsthm,latexsym}
\usepackage{amsfonts}
\usepackage{amsfonts}
\usepackage{graphicx}
\usepackage{textcomp}
\usepackage{cite}
\usepackage{enumerate}
\usepackage[mathscr]{euscript}
\usepackage{mathtools}
\newtheorem{theorem}{Theorem}[section]

\newtheorem{definition}[theorem]{Definition}

\newtheorem{problem}[theorem]{Problem}

\title{This is the title}
\usepackage{amssymb}
\usepackage{amssymb}
\usepackage{amssymb}
\usepackage{amssymb}
\usepackage{amsmath}
\usepackage{tikz}
\usepackage{hyperref}
\usepackage{enumerate}
\usepackage{mathtools}
\usepackage{amsmath}
\usepackage{tikz}
\usepackage{amssymb}
\usepackage{amsmath}
\usepackage{tikz}
\usepackage{hyperref}
\usepackage{amssymb}
\usepackage{amssymb}
\usepackage{amssymb}
\usepackage{amssymb}
\usepackage{amsmath}
\usepackage{tikz}
\usepackage{hyperref}
\usepackage{enumerate}
\usepackage{mathtools}
\usepackage{amsmath}
\usepackage{tikz}
\usepackage{graphicx}
\usepackage{textcomp}
\usetikzlibrary{cd}
\usepackage{fancyhdr}

\begin{document}
\begin{center}
{\bf{THE $(abc,pqr)$-PROBLEM FOR APPROXIMATE SCHAUDER FRAMES FOR  BANACH SPACES}}\\
\textbf{K. MAHESH KRISHNA} \\
Post Doctoral Fellow \\
Statistics and Mathematics Unit\\
Indian Statistical Institute, Bangalore Centre\\
Karnataka 560 059 India\\
Email: kmaheshak@gmail.com \\
Date: \today
\end{center}
\hrule
\vspace{0.5cm}
\textbf{Abstract}: Motivated from the complete solution of important $abc$-problem for Gabor system for the Hilbert space $\mathcal{L}^2(\mathbb{R})$ by Dai and Sun [\textit{Memoirs of Amer. Math. Soc., 2016}] and from the existential result of approximate Schauder frames for $\mathcal{L}^p(\mathbb{R})$ using translation operators on $\mathcal{L}^p(\mathbb{R})$ by Freeman, Odell, Schlumprecht, and Zsak [\textit{Israel. J. Math, 2014}], we formulate $(abc, pqr)$-problem for approximate Schauder frames for Banach spaces  $\mathcal{L}^p(\mathbb{R})$, $1<p<\infty$.

\textbf{Keywords}: $abc$-problem for Gabor system,  Approximate Schauder Frame.

\textbf{Mathematics Subject Classification (2020)}: 42C15, 46E30.

\tableofcontents 

\section{Introduction}
Let   $\mathcal{H}$ be a  Hilbert space over $\mathbb{C}$. Recall that \cite{DUFFIN, OLEBOOK, HANMEMOIRS}
a sequence $\{\tau_n\}_n$ in  $\mathcal{H}$ is said to be a
frame for $\mathcal{H}$ if there exist $a,b>0$ such that 
\begin{align*}
a\|h\|^2 \leq \sum_{n=1}^\infty |\langle h, \tau_n\rangle|^2\leq b\|h\|^2, \quad \forall h \in \mathcal{H}.
\end{align*}
Among various examples of frames for Hilbert spaces, one of the most studied classes of frames for the Hilbert space $\mathcal{L}^2(\mathbb{R})$ is the frames arising from the Gabor systems 
\begin{align*}
	\{E_{mb}T_{na}g\}_{m,n \in \mathbb{Z}},
\end{align*}
where $g\in \mathcal{L}^2(\mathbb{R})$ and for $c, d \in \mathbb{R}$,
\begin{align*}
&T_c:	\mathcal{L}^2(\mathbb{R})\ni f \mapsto T_cf \in \mathcal{L}^2(\mathbb{R});~  T_cf: \mathbb{R} \ni x \mapsto ( T_cf)(x)\coloneqq f(x-c)\in \mathbb{C},\\
&E_d:	\mathcal{L}^2(\mathbb{R})\ni f \mapsto E_df \in \mathcal{L}^2(\mathbb{R});~  E_df: \mathbb{R} \ni x \mapsto ( E_df)(x)\coloneqq e^{2\pi i d x}f(x)\in \mathbb{C}.
\end{align*}
A detailed analysis of Gabor frames for Hilbert spaces can be found in \cite{FLANDRIN, DAUBECHIESBOOK, FEICHTINGERSTROHMER, FEICHTINGERSTROHMER2, GOHRONSHEN, GROCHENIGMYSTERY, OKOUDJOU}.   A fundamental problem in Gabor analysis is to characterize parameters $a>0$, $b>0$ and the function $g$ such that $	\{E_{mb}T_{na}g\}_{m,n \in \mathbb{Z}}$ is a frame for $\mathcal{L}^2(\mathbb{R})$. The problem has been solved  fully or partially for certain classes of functions in \cite{JANSSEN1981, JANSSEN1982, FEICHTINGERKAIBLINGER, GROCHENIG2001, JANSSEN2000, DAUBECHIESGROSSMANNMEYER, GROCHENIGSTOCKLER, GUHAN, RONSHEN, JANSSENSTROHMER2002, JANSSEN1996, JANSSEN2003, LYUBARSKII, SEIPWALLSTEN1992, SEIP1992, LYUBARSKIISEIP, BORICHEVGROCHENIGLYUBARSKII, JANSSEN1994, DAUBECHIESGROSSMANN, GROCHENIGROMEROSTOCKLER, LYUBARSKIINES, LEMVIG2017, LEMVIGHAAHR, GUHAN, HELAU, CASAZZAKALTON, HEIL2007, JANSSEN1995, JANSSENAEJM, DAUBECHIES1990, LANDAU, KLOOSSTOCKLER}. Let $c>0$. When $g=\chi_{[0,c]}$, the characteristic function on $[0,c]$, the problem of characterizing  parameters $a>0$, $b>0$ and $c>0$ such that 
\begin{align*}
	\{E_{mb}T_{na}\chi_{[0,c]}\}_{m,n \in \mathbb{Z}}
\end{align*}
is a frame for $\mathcal{L}^2(\mathbb{R})$ is known as $abc$-problem for Gabor system. A list of partial answers for this are obtained in \cite{LEMVIGHAAHR, GROCHENIGMYSTERY, LEMVIG2017, LYUBARSKIINES, GUHAN, HELAU, HEIL2007, CASAZZAKALTON, FEICHTINGERSTROHMER2}. In 2013, Dai and Sun fully resolved this problem in their Memoirs  
\cite{DAISUN, DAISUNBOOK}. In this paper, we formulate  a similar problem for approximate Schauder frames for function spaces.

\section{The $(abc, pqr)$-Problem for Approximate Schauder Frames}
We are motivated from the following two results.
\begin{enumerate}[\upshape(i)]
	\item The complete solution of $abc$-problem for Gabor system by Dai and Sun \cite{DAISUN}.
	\item For every $2<p<\infty$, there is an  approximate Schauder frame for $\mathcal{L}^p(\mathbb{R})$ obtained  using translation operators on $\mathcal{L}^p(\mathbb{R})$ \cite{FREEMANODELL}.
\end{enumerate}
  Generalized Fourier expansion which resulted from the theory of frames for Hilbert spaces led the notion of framing in \cite{CASAZZAHANLARSON} which is generalized to Schauder frames in \cite{CASAZZA}   which is further generalized in \cite{FREEMANODELL, THOMAS} to approximate Schauder frames defined as follows.
\begin{definition}\cite{FREEMANODELL, THOMAS}\label{ASFDEF}
 Let $\mathcal{X}$ be a separable Banach space and $\mathcal{X}^*$ be its dual.	Let $\{\tau_n\}_n$ be a sequence in  $\mathcal{X}$ and 	$\{f_n\}_n$ be a sequence in  $\mathcal{X}^*.$ The pair $ (\{f_n \}_{n}, \{\tau_n \}_{n}) $ is said to be an \textbf{approximate Schauder frame (ASF)} for $\mathcal{X}$ if the map 
	\begin{align*}
	S_{f, \tau}:\mathcal{X}\ni x \mapsto S_{f, \tau}x\coloneqq \sum_{n=1}^\infty
	f_n(x)\tau_n \in
	\mathcal{X}
	\end{align*}
	is a well-defined bounded linear, invertible operator.
\end{definition} 
We now consider the Banach space $\mathcal{X}=\mathcal{L}^p(\mathbb{R})$, for $1<p<\infty$. Let $q$ be the conjugate index of $p$. Since dual of $\mathcal{L}^p(\mathbb{R})$ is isometrically isomorphic to $\mathcal{L}^q(\mathbb{R})$, given $\phi \in \mathcal{X}^*$, there exists a unique $\omega\in \mathcal{L}^q(\mathbb{R})$ such that 
\begin{align*}
\phi(u)=\int_\mathbb{R}u(\alpha)\omega(\alpha)\,d\alpha \eqqcolon [u,\omega], \quad \forall u \in \mathcal{L}^p(\mathbb{R}).
\end{align*}

  In particular, if $\{f_n\}_n$ is a sequence in  $\mathcal{X}^*$, then there exists a unique sequence $\{\omega_n\}_n$ in  $\mathcal{X}$ such that 
  \begin{align*}
  	f_n(x)=[x, \omega_n], \quad \forall x \in \mathcal{X}, \forall n \in \mathbb{N}.
  \end{align*}
 Thus asking whether a pair $ (\{f_n \}_{n}, \{\tau_n \}_{n}) $ is  an ASF for $\mathcal{X}$ is same as asking whether the map 
 	\begin{align*}
 	S_{\omega, \tau}:\mathcal{X}\ni x \mapsto S_{\omega, \tau}x\coloneqq \sum_{n=1}^\infty
 [x, \omega_n]	\tau_n \in
 	\mathcal{X}
 \end{align*}
 is a well-defined bounded linear, invertible operator. We then have the following problem.
   \begin{problem}\label{GABORBANACH}
 \textbf{Let $1<p<\infty$ and $q$ be the conjugate index of $p$. Characterize parameters $a>0$, 	$b>0$, $p>0$, $q>0$ and functions $g\in \mathcal{L}^p(\mathbb{R})$, $h\in \mathcal{L}^q(\mathbb{R})$ such that the pair 
  \begin{align*}
  	(\{\tilde{E}_{mq}\tilde{T}_{np}h\}_{m,n \in \mathbb{Z}}, 	\{E_{mb}T_{na}g\}_{m,n \in \mathbb{Z}}),
  \end{align*}
is an ASF for $\mathcal{L}^p(\mathbb{R})$, where 
\begin{align*}
	&T_c:	\mathcal{L}^p(\mathbb{R})\ni f \mapsto T_cf \in \mathcal{L}^p(\mathbb{R});~  T_cf: \mathbb{R} \ni x \mapsto ( T_cf)(x)\coloneqq f(x-c)\in \mathbb{C}, ~ \text{ for } c\in \mathbb{R},\\
	&E_d:	\mathcal{L}^p(\mathbb{R})\ni f \mapsto E_df \in \mathcal{L}^p(\mathbb{R});~  E_df: \mathbb{R} \ni x \mapsto ( E_df)(x)\coloneqq e^{2\pi i d x}f(x)\in \mathbb{C}, ~ \text{ for } d\in \mathbb{R}
\end{align*}
and 
\begin{align*}
	&\tilde{T}_s:	\mathcal{L}^q(\mathbb{R})\ni f \mapsto \tilde{T}_sf \in \mathcal{L}^q(\mathbb{R});~  \tilde{T}_sf: \mathbb{R} \ni x \mapsto ( \tilde{T}_sf)(x)\coloneqq f(x-s)\in \mathbb{C}, ~ \text{ for } s\in \mathbb{R},\\
	&\tilde{E}_t:	\mathcal{L}^q(\mathbb{R})\ni f \mapsto \tilde{E}_tf \in \mathcal{L}^q(\mathbb{R});~  \tilde{E}_tf: \mathbb{R} \ni x \mapsto ( \tilde{E}_tf)(x)\coloneqq e^{2\pi i t x}f(x)\in \mathbb{C}, ~ \text{ for } t\in \mathbb{R}.
\end{align*}}
  \end{problem}
Note that all operators $T_c$, $E_c$,  $\tilde{T}_s $,  $\tilde{E}_t $ are invertible isometries. We also note that Problem \ref{GABORBANACH} can be stated for space  $\mathcal{L}^p[a,b].$ A particular case of Problem  \ref{GABORBANACH}  is the following $(abc, pqr)$-problem.
  \begin{problem}\textbf{(The $(abc, pqr)$-problem for approximate Schauder frames)}\label{ABCPQR}
   \textbf{Let $1<p<\infty$ and $q$ be the conjugate index of $p$. Characterize parameters $a>0$, 	$b>0$, $c>0$, $p>0$, $q>0$, $r>0$   such that the pair 
  	\begin{align*}
  		(\{\tilde{E}_{mq}\tilde{T}_{np}\chi_{[0,r]}\}_{m,n \in \mathbb{Z}}, 	\{E_{mb}T_{na}\chi_{[0,c]}\}_{m,n \in \mathbb{Z}}),
  	\end{align*}
  	is an ASF for $\mathcal{L}^p(\mathbb{R})$}.	
  \end{problem}
 We remark that we can formulate problems similar to Problems \ref{GABORBANACH}, \ref{ABCPQR}  for p-approximate Schauder frames, p-approximate Bessel sequences,  p-approximate orthonormal bases and p-approximate Riesz bases \cite{KRISHNAJOHNSON, KRISHNA, KRISHNAJOHNSON2, KRISHNAJOHNSON3}. Following is one such in the discrete case. Our problem is motivated from the discrete Gabor analysis in $\ell^2(\mathbb{Z})$ \cite{LOPEZHAN, CVETKOVICVETTERLI, MORRISLU, BOLCSKEIHLAWATSCH, CVETKOVICVETTERLI2, BANNERTGROCHENIGSTOCKLER}.
 \begin{problem}\textbf{(The $(MN, PQ)$-problem for p-approximate Schauder frames)}
 \textbf{Let $1<p<\infty$ and $q$ be the conjugate index of $p$. Characterize natural numbers $N,M,P,Q$ and sequences $\{x_j\}_{j\in \mathbb{Z}}\in \ell^p(\mathbb{Z})$, $\{y_j\}_{j\in \mathbb{Z}}\in \ell^q(\mathbb{Z})$ such that the pair 
 	\begin{align*}
 	(\{	\tilde{E}_{m/P}\tilde{T}_{nQ}\{y_j\}_{j\in \mathbb{Z}}\}_{0\leq m\leq P-1, n\in \mathbb{Z}}, \{E_{m/M}T_{nN}\{x_j\}_{j\in \mathbb{Z}}\}_{0\leq m\leq M-1, n\in \mathbb{Z}})
 	\end{align*}
 is a 
 \begin{enumerate}[\upshape(i)]
 	\item p-orthonormal basis for $ \ell^p(\mathbb{Z})$.
 	\item p-approximate Riesz basis for $ \ell^p(\mathbb{Z})$.
 	\item p-approximate Schauder frame for $ \ell^p(\mathbb{Z})$.
 	\item p-approximate Bessel sequence for $ \ell^p(\mathbb{Z})$, 
 \end{enumerate}
where
\begin{align*}
&	\tilde{E}_{m/P}\tilde{T}_{nQ}:  \ell^q(\mathbb{Z})\ni \{a_j\}_{j\in \mathbb{Z}}\mapsto \tilde{E}_{m/P}\tilde{T}_{nQ}\{a_j\}_{j\in \mathbb{Z}}\coloneqq  \{e^{2\pi ijm/P}a_{j-nQ}\}_{j\in \mathbb{Z}}\in  \ell^q(\mathbb{Z}),\\
&	E_{m/M}T_{nN}:  \ell^p(\mathbb{Z})\ni \{b_j\}_{j\in \mathbb{Z}}\mapsto E_{m/M}T_{nN}\{a_j\}_{j\in \mathbb{Z}}\coloneqq  \{e^{2\pi ijm/M}b_{j-nN}\}_{j\in \mathbb{Z}}\in  \ell^p(\mathbb{Z}).
\end{align*}}
 \end{problem} 
 Next we formulate the continuous version of Problem  \ref{GABORBANACH} (see  \cite{KRISHNA2} for the definition of continuous approximate Schauder frame and \cite{EISNERFREEMAN} for the definition of continuous Schauder frame).
  \begin{problem}
   \textbf{Let $1<p<\infty$ and $q$ be the conjugate index of $p$. Characterize  functions $g\in \mathcal{L}^p(\mathbb{R})$, $h\in \mathcal{L}^q(\mathbb{R})$ such that the pair 
  	\begin{align*}
  		(\{	\tilde{E}_b	\tilde{T}_ah\}_{a,b \in \mathbb{R}}, \{	E_b	T_ag\}_{a,b \in \mathbb{R}})
  	\end{align*}
  is a continuous approximate Schauder frame for $\mathcal{L}^p(\mathbb{R})$}.
\end{problem}

 \bibliographystyle{plain}
 \bibliography{reference.bib}

\end{document}